\documentclass[12pt]{article}
\usepackage[english,russian]{babel}
\usepackage{latexsym, amsgen, amsmath, amstext, amsbsy, amsopn, amsfonts, amsthm, amssymb}
\usepackage[symbol]{footmisc}

\renewenvironment{abstract}{\vspace{-1cm}\small\quotation}{\endquotation}
 \RequirePackage[all]{xy}

\newenvironment{network}%
    {\catcode`"=12\begin{xy}<1ex,0ex>:}%
    {\end{xy}\catcode`"=13}
\def\nnNode"#1"(#2)#3{\POS(#2)*#3="#1"}
\def\nnLink"#1,#2"#3{\POS"#1"\ar #3 "#2"}

\begin{document}
\newpage
\centerline{\textbf{Constructive Graph Theory:}}
\centerline{\textbf{Generation Methods, Structure  }}
\centerline{\textbf{and Dynamic Characterization }}
 \centerline{\textbf{of Closed Classes
of Graphs - A survey}}

\medskip
\centerline { Mikhail Iordanski \footnote[1] {Lobachevsky State
University, Minin Pedagogical State University of Nizhny Novgorod,
Russia, iordanski@mail.ru}}
 \bigskip
\bigskip
\begin{abstract}
\centerline { Abstract}
 The processes of constructing some graphs
from others using binary operations of union with intersection
(gluing) are studied. For graph classes closed with respect to
gluing operations the elemental and operational bases are
introduced. The genera\-ting bases together with the system of
restrictions on the gluing operation, that preserve the
characteristic properties of
 graphs
form a constructive descriptions of the closed classes of
 graphs.
 It is shown that each
closed class of graphs has a unique elemental basis and at least
one operational basis. For the closed class of all graphs and all
basis precomplete closed subclasses of it the constructive
descriptions are considered. For each of them its characteristic
properties and a diagram of the inclusion of subclasses in
superclasses is given. Constructive descriptions are obtained for
some classes of graphs with classical properties. Some possible
applications of constructive theory are discussed in conclusion
\end{abstract}

\leftline{{\bf{\large 1. Introduction}}}
\medskip

 The work is a systematic review of the results on
the constructive theory of graphs from Russian-language
publications, mainly by the author, over the last 25 years.

The most complete description of the results is given in the
monograph [1]. In the monograph and in the articles
 [2-5] different
systems of operations on graphs are considered. In this work we
based on the using binary set-theoretical operations of union with
intersection (gluing operations), which is the most natural
representation.

In the constructive theory, the processes of building some graphs
from others are being studied. The generation methods, structure
and dynamic characterizat\-ion of graph classes closed with
respect
 to gluing operations are considered.
Restrictions on gluing operations
 are studied, under which various characteristic
properties of graphs are saved. Knowing such restrictions allows
you to build graphs that have the specified properties.

In the general case, loops and multiple edges are permissible in
graphs. The following notation is used: $K_{n}$ -- complete graph,
$C_{n}$ -- simple cycle, $L_{n}$ -- simple chain, $O_{n}$ -- empty
graph, all of them contain $n$ vertices ($O_{o}$ -- null graph
without vertices).

Let $G_{1}$ and $G_{2}$ are disjoint graphs.
The \emph{gluing operation} consists in the identification
of
isomorphic subgraphs $G_{1}' \subseteq~G_{1} $ and $G_{2}'
\subseteq G_{2}$. Glueing operation is called \emph{trivial}
 if
$G_{1}'= G_{1}$ and (or) $G_{2}' = G_{2}$. For each of the graphs
obtained as a result of the operation  of gluing the graphs
$G_{1}$ and $G_{2}$ on the subgraph $\tilde{G}$, isomorphic to
$G_{1}'$ and $G_{2}'$ the notation $(G_{1}\,\circ
\,G_{2})\,\tilde{G}$ is used. Operand graphs $G_{1}$ and $G_{2}$
are isomorphic to subgraphs of the resulting graph $G =
(G_{1}\,\circ \,G_{2})\,\tilde{G}$ of
 gluing
operation. The subgraph $\tilde{G}$ is called the
\emph{subgraph of
gluing}. For fixed graphs $G_{1}$ and $G_{2}$ the resulting
graph
$G = (G_{1}\,\circ \,G_{2})\,\tilde{G}$ may depend on the
type of
subgraph gluing $\tilde{G}$ (Fig 1), the choice of identifiable
subgraphs  $G_{1}'$ and $G_{2}'$ in operand graphs (Fig 2) and the
method of their identification (Fig 3).

\begin{center}

\unitlength=1.4mm \special{em:linewidth 0.4pt}
\linethickness{0.4pt}

\begin{picture}(100,27)
\put(35,3){\circle*{0.8}} \put(30,13){\circle*{0.8}}
\put(25,3){\line(1,0){10}} \put(25,3){\line(1,2){5}}
\put(35,3){\line(-1,2){5}}
\put(10,7){\makebox(0,0)[cb]{$G_2=K_3$}}
\put(25,18){\circle*{0.8}} \put(35,18){\circle*{0.8}}
\put(30,28){\circle*{0.8}} \put(25,18){\line(1,0){10}}
\put(25,18){\line(1,2){5}} \put(35,18){\line(-1,2){5}}
\put(10,22){\makebox(0,0)[cb]{$G_1=K_3$}}
\put(65,17){\circle*{0.8}} \put(65,27){\circle*{0.8}}
\put(85,27){\circle*{0.8}} \put(85,17){\circle*{0.8}}
\put(75,22){\circle*{0.8}} \put(65,17){\line(0,1){10}}
\put(65,17){\line(2,1){20}} \put(85,17){\line(-2,1){20}}
\put(85,17){\line(0,1){10}} \put(50,14){\vector(2,-3){7}}
\put(50,16){\vector(1,1){7}} \put(43,10){\oval(4,10)[br]}
\put(47,10){\oval(4,10)[tl]} \put(43,20){\oval(4,10)[tr]}
\put(47,20){\oval(4,10)[bl]} \put(65,5){\circle*{0.8}}
\put(85,5){\circle*{0.8}} \put(75,0){\circle*{0.8}}
\put(75,10){\circle*{0.8}} \put(65,5){\line(2,1){10}}
\put(65,5){\line(2,-1){10}} \put(85,5){\line(-2,1){10}}
\put(85,5){\line(-2,-1){10}} \put(75,0){\line(0,1){10}}
\put(59,16){\makebox(0,0)[cb]{$\tilde{G}=O_1$}}
\put(59,10){\makebox(0,0)[cb]{$\tilde{G}=K_2$}}
\put(55,24){\makebox(0,0)[cb]{$(K_3\,\circ\,K_3)\,O_1$}}
\put(55,0){\makebox(0,0)[cb]{$(K_3\,\circ\,K_3)\,K_2$}}
\put(25,3){\circle*{0.8}}
\end{picture}
\end{center}
\medskip
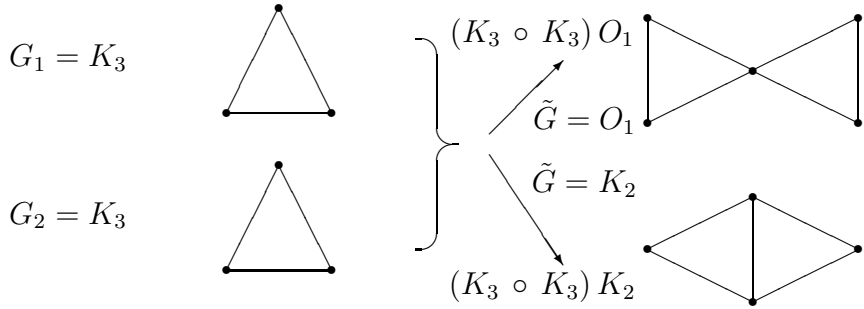
Figure 1. The result of the gluing operation depends on the type
of the gluing subgraph~$\tilde{G}$.

\begin{center}
\unitlength=1.4mm \special{em:linewidth 0.4pt}
\linethickness{0.4pt}
\begin{picture}(100,29)
\put(43,9){\oval(4,8)[br]}
\put(58,12){\makebox(0,0)[cb]{$\tilde{G}= O_2$}}
\put(74,13){\makebox(0,0)[cb]{$(L_3\,\circ\,K_3)\,O_2$}}
\put(50,12){\vector(1,-1){7}} \put(50,14){\vector(2,3){7}}
\put(47,9){\oval(4,8)[tl]} \put(43,17){\oval(4,8)[tr]}
\put(47,17){\oval(4,8)[bl]}
\put(10,20){\makebox(0,0)[cb]{$G_1=L_3$}}
\put(10,9){\makebox(0,0)[cb]{$G_2=K_3$}} \put(25,5){\circle*{0.8}}
\put(35,5){\circle*{0.8}} \put(30,15){\circle*{0.8}}
\put(20,20){\circle*{0.8}} \put(30,20){\circle*{0.8}}
\put(40,20){\circle*{0.8}} \put(20,20){\line(1,0){20}}
\put(25,5){\line(1,0){10}} \put(25,5){\line(1,2){5}}
\put(35,5){\line(-1,2){5}} \put(65,0){\circle*{0.8}}
\put(75,0){\circle*{0.8}} \put(85,0){\circle*{0.8}}
\put(70,10){\circle*{0.8}} \put(65,0){\line(1,0){20}}
\put(65,0){\line(1,2){5}} \put(75,0){\line(-1,2){5}}
\bezier{125}(70,10)(76,8)(75,0) \put(62,25){\circle*{0.8}}
\put(82,25){\circle*{0.8}} \put(72,20){\circle*{0.8}}
\put(72,30){\circle*{0.8}} \put(62,25){\line(2,1){10}}
\put(62,25){\line(2,-1){10}} \put(82,25){\line(-2,1){10}}
\put(82,25){\line(-2,-1){10}} \put(72,20){\line(0,1){10}}
\end{picture}
\end{center}

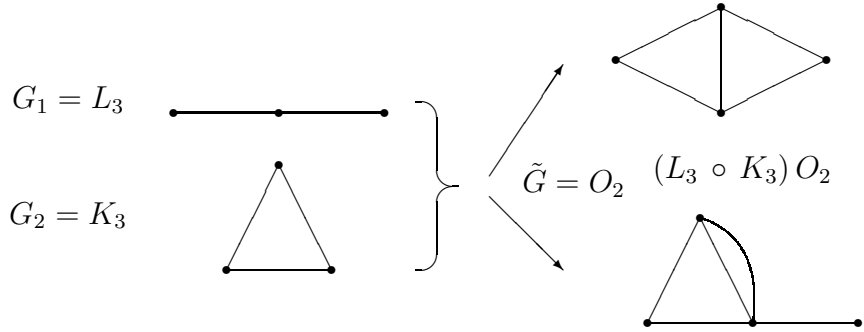
Figure 2. The result of the gluing operation on subgraph
$\tilde{G}= O_2$ depends on the choice of identifiable vertices in
operand graphs $G_{1}$ and $G_{2}$.

\begin{center}
\unitlength=1.4mm \special{em:linewidth 0.4pt}
\linethickness{0.4pt}
\begin{picture}(100,27)
\put(50,18){\vector(1,1){7}} \put(50,16){\vector(3,-1){7}}
\put(43,13.5){\oval(4,7)[br]} \put(47,13.5){\oval(4,7)[tl]}
\put(43,20.5){\oval(4,7)[tr]} \put(47,20.5){\oval(4,7)[bl]}
\put(10,22){\makebox(0,0)[cb]{$G_1=L_3$}}
\put(10,12){\makebox(0,0)[cb]{$G_2=L_3$}}
\put(20,12){\circle*{0.8}} \put(30,12){\circle*{0.8}}
\put(40,12){\circle*{0.8}} \put(20,22){\circle*{0.8}}
\put(30,22){\circle*{0.8}} \put(40,22){\circle*{0.8}}
\put(20,12){\line(1,0){20}} \put(20,22){\line(1,0){20}}
\put(65,25){\line(1,0){30}} \put(65,25){\circle*{0.8}}
\put(75,25){\circle*{0.8}} \put(85,25){\circle*{0.8}}
\put(95,25){\circle*{0.8}}

\put(70,0){\circle*{0.8}} \put(80,0){\circle*{0.8}}
\put(75,5){\circle*{0.8}} \put(75,12){\circle*{0.8}}
\put(75,5){\line(0,1){7}} \put(70,0){\line(1,1){5}}
\put(80,0){\line(-1,1){5}}

\put(77,16){\makebox(0,0)[cb]{$(L_3\,\circ\,L_3)\,K_2$}}

\put(57,16){\makebox(0,0)[cb]{$\tilde{G}=K_2$}}
\end{picture}
\end{center}
\medskip

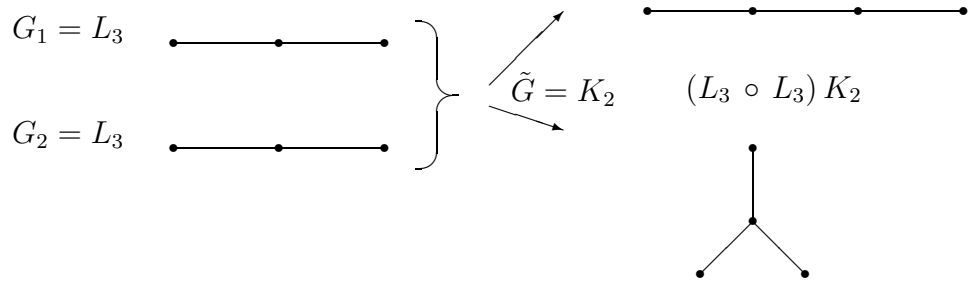
 Figure 3. The result of the gluing operation on subgraph $\tilde{G}= K_2$ with fixed subgraphs $K_2$
 in operand graphs depends on the choice of  the method of their identification.

\medskip

Let $\Im$ be a set of graphs. Graph $G$ is called {\it
superposition} of graphs from $\Im$ if $G \in \Im$ or $G$ can
 be
obtained by successive application of the operations of gluing
 to
the graphs from $\Im$ or to the graphs obtained from $\Im$
with
the operations of gluing. When performing each gluing
operation,
the type of the identified subgraphs their choice in operand
graphs
  and the method of identification are
determined independently. The process of construc\-tion graph
 $G$
from the graphs of the set $\Im$ determines
 {\it the superposition
operation} of graphs from $\Im$. If in the superposition
 operation
at least one of the graph-operands of each gluing operation
belongs to the set of $\Im$, then the superposition operation
is called \emph{canonical}.

The set of  graphs $\Im $, as well as graphs, derived from
 $ \Im $
using superposition operations, denoted by $ [\Im] $. If $
[\Im] =
\Im $, then the class $ \Im $ is called \emph {closed}. A
closed
class of graphs is called \emph {trivial} if used in
superposition
operations only trivial gluing operations. The set of
 $ \Im '
\subset \Im $ is \emph {a complete system of graphs} in a
closed
class $ \Im $ if $ [\Im '] = ~ \Im $. Minimal on inclusion
the
complete system of graphs  $ B_e $ is called
 \emph {elemental
basis} of a closed class $ \Im $.

 Operations with isomorphic gluing subgraphs $ \tilde {G} $
refer to one \emph {type}. The set of types gluing operations,
use
of which is enough to build from $ B_e $ all graphs of a
closed
class $ \Im $, forms \emph {a complete  system of types
 gluing
operations}. Minimal on inclusion the complete  system of
types
gluing operations $ B_o $ is called \emph {operational basis}
 of
closed class $ \Im $. The operational basis $ B_o $ is
described
by set of graphs isomorphic to subgroups of gluing
$ \tilde {G} $.
The elemental and operational bases are called
\emph {generators
bases}. The generating bases define \emph {constructive
description} of the closed class $ \Im $.

The numbering of all statements given in the text
(theorems,
lemmas and corollaries) is independent in each subsection.
The section number is indicated first, then the subsecti\-on
number.

Proofs of the fundamental statements of the first two
subsections of the the main part is given.
Ends of the proofs are marked with $\Box$.

\medskip
\medskip
  \leftline{{\bf {\large   2. The generation methods,
 structure and dinamic  }}}
 \leftline{{\bf {\large \hspace {0.5 cm} characterizations of
 closed classes of graphs
 }}}
\medskip

\leftline{{\bf 2.1. Generation Methods of closed classes of
 graphs}}
\medskip
\medskip

A number of theorems on the structure and methods of
generating closed classes of graphs were announced without
proofs in [6]. In subsequent works, all of them were proved.
 Some of the most important ones are listed below.

\medskip
{\bf Theorem 2.1.1} [7].
{\it Every closed class of graphs
$ \Im $ has single elemental basis.}
\medskip

 {\bf Proof.} Consider an arbitrary closed class of graphs $ \Im $.
We associate with it an infinite oriented graph $ G ^ \Im $.
Its vertices correspond to the graphs from $ \Im $. Arc
$ (v_i,v_j) $,
where $ v_i, v_j \in V (G ^ \Im) $, $ i \neq j $, is carried
out
then and only when the graph $ G_i $ corresponding to the
vertex $
v_i $ is graph operand of at least one non-trivial operation
binary gluing that implements the graph $ G_j $,
corresponding to vertex $ v_j $.
Since the arcs of the graph $ G ^ \Im$
correspond only to non-trivial gluing operations, then
$|V(G_j)|>|V(G_i)|$ or (and) $ | E (G_j) |> | E (G_i) |$.
From the finiteness of each graph in $ \Im $ it follows that
all paths leading
to any vertex of the graph ~ $ G ^ \Im $, contain a finite
number of different vertices.
Graph $ G ^ \Im $ cannot have oriented simple cycles since
their vertices would correspond to isomorphic graphs.
 So
all the paths leading to any vertex the graph ~ $ G ^ \Im $,
 have
finite length. It follows that the set vertices of the graph
$ G ^
\Im $ with indegrees equal to zero is not empty.
The graphs corresponding to such vertices form the elemental
basis
of a closed class of graphs ~ $ \Im $, since none of these
graphs
can be expressed as \hbox {superposition} of other graphs
from $
\Im $. The uniqueness of the elemental basis follows
 from its definition. $\Box$

\medskip
{\bf Theorem 2.1.2} [7].
 {\it Power of the set of all closed classes graphs
are continual.}
\medskip

{\bf Proof.} The number of closed classes of graphs can be
estimated from above the number of all subsets of the countable
set  of all graphs. For lower estimates it is enough to select an
infinite sequence graphs, each of which cannot be represented by a
superposition of others sequence graphs, for example, $ C_n, n =
1,2, ... \,.$ Choosing all possible subsets of this sequence as
elementary bases of the corresponding closed classes, we obtain a
continuous set of closed classes.

\medskip
{\bf Corollary 2.1.1} {\it There are closed classes of
 graphs
with $ |B_e | = \infty $.}

For operational bases, was obtained the following result.

\medskip
{\bf Theorem 2.1.3} [8].
 {\it Every closed class of graphs $ \Im $
 has an operational basis.}
 \medskip

{\bf Proof.} From the definition of a complete system of types of
gluing operations it follows that each closed class of graphs
$\Im$ has a nonempty set  of such systems. This set contains, for
example, a system including subgraphs of all graphs from $ \Im $.
If graphs from $\Im $ can be constructed using canonical
 superpositions then we have also the complete system
including only subgraphs of graphs from $B_e$.
Full can also be subgraph systems of various other subsets of
graphs from $ \Im $.

Each complete system of types of gluing operations is defined
by the set of graphs.
Put each graph $ G $ in one-to-one match the
positive
integer $ n (G) $ so that no graph with a higher number would
 not
be isomorphic to a subgraph of a graph with a smaller number.
 This
can always be done, for example, by numbering the graphs in
the
non-decreasing order of the sum of the number of their
 vertices
and edges. Graph, corresponding to the number $n$, we denote
by $ G (n) $.

 We can assign the characteristic binary
fraction $ 0, r_ {1} r_ {2} ... r_ {n} ... $ to each set of graphs
$ R $, in which $r_ {n} = 0 $ if $ G (n) $ does not belong to the
set $ R $ and $r_ {n} = 1 $ if $ G (n) $ belongs to the set $ R $.
According to this rule we can associate binary fraction as real
number with each complete system of types of gluing operations.

The set of all complete systems of types of gluing operations
corresponds to the set $A$ of real numbers.
Since these numbers are positive, there exists a number
$\inf A $.
Show that the number $\inf A $ also belongs to the set A,
that is, it corresponds to the complete system of types of
gluing operations.
Suppose that the system corresponding to the
number $ \inf A $ is not complete. Then there is a graph
$ G \in
\Im $ that cannot be constructed from graphs of elemental
 basis
$B_e $ using a system of types of gluing operations
corresponding to the number $ \inf A $.

From the definition of infinium, it follows that
there is a complete system of types of gluing
operations, which corresponds to the
number $ M < (\inf A + 2 ^ {- n (G)}) $. Since the graphs
with
numbers greater than $ n (G) $ are not can be used in the
construction of the graph $ G $ (they are not isomorphic to
own
subgraphs of $ G $), and the system corresponding to the
number
$ M $ is
complete, then the graph $ G $ can be constructed using
system of types of gluing operations
corresponding to the number of $ \inf A $.

A complete system of types of gluing operations corresponding
to
the number $ \inf A $ is minimal on inclusion, since any of
its
own subset corresponds to a smaller number, but the
complete
 system of types of gluing operations corresponding to a
number less than $ \inf
A $ does not exist. Thus, the system of types of operations
corresponding to the number of $ \inf A $, is the operational
basis of the closed class of graphs $ \Im $. $\Box$

The resulting graph of any gluing operation saves such
properties
of graphs-operands as the absence of isolated vertices,
loops or edges.

It is not difficult to see that if all the graphs from
elemental basis $ B_e $ are connected, then to obtain
disconnected
graphs it is necessary to include a null graph $ O_0 $ in the
operational basis $ B_o $.

\medskip
{\bf Theorem 2.1.4} [7]. {\it Closed class of all graphs  has
an elemental basis
 $ B_ {e} = \{O_1, C_1, K_2 \} $ and operational
basis $ B_ {o} = \{O_0, O_1, O_2 \} $.}
 \medskip

 {\bf Proof.} As each gluing operation saves in graphs the lack of
   isolated vertices, loops and edges, we have the inclusion
$ \{O_1, C_1, K_2 \} \subset B_e $. Since the elemental basis
of
each a closed class of graphs is unique, then for proof
reverse
inclusion enough to show that any graph $ G $ we represent as
a
superposition of graphs from set $ \{O_1, C_1, K_2 \} $. It
can be
done, for example, like this:

1) construct an empty $ | V (G) | $-vertex graph using $(|V(G)|-1)
$ gluing operati\-ons on $ O_0 $ implementing graphs of the form $
(g \circ O_1) O_0 $, where $ g $ is the resulting graph of the
previous one gluing operations ($ g = O_1 $ when performing the
first operation);

2) supplement the empty graph with edges up to the graph $G$,
using $|E(G)|$ gluing operations that implement graphs of the
form $ (g \circ K_2) O_2 $ and (or) graphs of the form
$ (g \circ C_1) O_1 $.

To complete the proof, it suffices to establish minimal on
inclusion the numbers of graphs included in the set $ B_o $.
 Without gluing operations on $ O_0 $ can not be implemented
disconnected graphs since all graphs of elemental basis  are
connected. Graphs containing multiple loops and edges cannot be
constructed without using gluing operations on $ O_1 $ and
$ O_2 $ respectively. $\Box$

In the proof of Theorem 2.1.4, only canonical superpositions
were used. This way of constructing graphs is always
admissible in the following case.

\medskip
{\bf Lemma 2.1.1} [9].
 {\it Graphs of a closed class $ \Im $ with
generators bases $ B_e $ and $ B_ {o} $ can be constructed
using
canonical superpositions if the operational basis
 $ B_ {o} =
\{O_0, O_1, \dots, O_n \} $, where $ n = \max_ {G \in B_e} | V(G)
| $.} \smallskip

{\bf Proof.} Each superposition operation that implements
arbitrary graph from $ \Im $,
 you can match its coverage to the $ B_e $ graphs. Consider
the graph
 cover. Its vertices are subgraphs isomorphic to the graphs
from $ B_e $ and edges join
 vertices corresponding to intersecting subgraphs.

 In a connected graph of a coverage containing at least two
vertices, always
 there is a vertex whose removal preserves connectivity. Any
vertex deletion process
 from the graph of the covering preserving the connectivity,
 the operation of canonical
 superposition corresponds to its reverse consideration.

Since when using gluing on empty subgraphs
 intersect, then for any order of graph assembling, all gluing
 subgraphs will be empty.
 Only the number of identified vertices in each specific
operation can vary.
 When using canonical superpositions, it cannot exceed the
number
 of  vertices in the added graph $ G \in B_e $, therefore,
 all operations satisfy the conditions of the lemma.

For disconnected coverage graphs, firstly using canonical
superposition of gluing operations
 a graph is constructed on $ O_0 $, each connected component
 of which is
isomorphic to the graph $ G \in  B_e $, which is the original in
the canonical superposition realizing this component. $\Box$

A canonical superposition is always possible if the gluing
operation has the property of associativity.
Restrictions under which the gluing operation is associative
were considered in [10]. There, in particular, the
associativity of operations over complete subgraphs of
gluing was shown.

\medskip

\leftline {{\bf 2.2. Structure of closed classes of graphs}}

\medskip

In the theory of functional systems with operations the
concept of
studying the structure of closed classes of boolean functions
 by using
precomplete classes is considered [11].
Closed class $ \Im_1 \subset \Im_2 $
called \emph {precomplete} in a closed superclass $ \Im_2 $
if $[\Im_1] \neq \Im_2 $, but adding to $ \Im_1 $ of any
element $ r\in \Im_2 \setminus \Im_1 $ we get $ [\Im_1 \cup r]
 = \Im_2 $. The precomplete class $ \Im_1 $ is called
\emph {trivial} if the set $\Im_2 \setminus \Im_1 $ contains
exactly one element. For closed classes of graphs the concept of a
precomplete class is not informative.

 \medskip
{\bf Theorem 2.2.1} [5].
{\it All precomplete closed classes of
graphs are trivial.}
 \medskip

{\bf Proof.} Assume that the subclass $ \Im_1 $ does not contain
two the graphs $ G_1 $ and $ G_2 $ from the superclass $ \Im_2 $.
If $ | V (G_ {1}) | <| V (G_ {2}) | $ and (or) $ | E (G_ {1}) | <|
E (G_ {2}) | $, then the graph $ G_ {1} $ cannot be built using
the graph $ G_ {2} $ because the graphs are operands gluing
operations are isomorphic to subgraphs of the resulting graph.
 If $ | V(G_ {1}) | = | V (G_ {2}) | $ and
$ | E (G_ {1}) | = | E (G_ {2})| $, but $ G_ {1}
 \ncong G_ {2}$, then by the same reason, none
of the graphs $ G_1 $ and $ G_2 $ cannot be constructed using
another graph.
Thus, all precomplete closed classes of graphs do not contain
 only one graph from their superclasses and are trivial.
$\Box$

 To describe the structure of closed classes of graphs, we
introduce the concept of basis precompleteness [4]. Class $ \Im_1
$ is \emph{precomplete in elemental basis} in $ \Im_2 $ if $ B_e $
of the class $ \Im_1 $ does not contain only one of the graphs of
elemental basis of the class $ \Im_2 $ and the operational bases
of both classes coincide. Similarly, the class $ \Im_1 $ is \emph
{ precomplete in operational basis} in $ \Im_2 $ if $ B_o $ of the
class $ \Im_1 $ does not contain only one of the graphs of the
operational basis of the class $ \Im_2 $ and the elemental bases
of both classes coincide.
\medskip

{\bf Theorem 2.2.2} [9].
{\it The closed class of all graphs
contains 35
nontrivial closed subclasses that are bases precomplete on
elemental or operational basis in their superclasses.}
\medskip

{\bf Proof.} Consider the closed classes of graphs whose
generators bases are subsets of the bases $ B_e $ or(and) $ B_o $
of closed class of all graphs. Constructive descriptions and the
characteristic properties of these classes are listed in the table
1 for connected graphs and in table 2 for graphs that admit
different number of connected components.

Number of vertices and edges in graphs (subgraphs) are denoted as
$ N (n) $ and $ M (m) $, respectively. Characteristic
properties
are formulated on analysis of the generating bases used in
the construction of graphs.

The union of the graphs $ G_1 $ and $ G_2 $ without intersections
obtained with the help of gluing operations $ (G_1 \, \circ \,
G_2) \, O_0 $. Adding edges and loops to the current graph $ G$ is
implemented, respectively, by gluing operations $ (G \, \circ \,
K_2) \, O_2 $ and $ (G \, \circ \, C_1) \, O_1 $. Adding edge with
vertice to the current graph $ G $ can be done by one of the
following gluing operations $ (G \, \circ \, K_2) \, O_1 $ or
$(((G \, \circ \, O_1) \, O_0, \circ \, K_2) \, O_2 $. Considering
Lemma 2.1.1, the construc\-ting of graphs can be restricted to
canonical superpositions. As elemental bases various subsets of
graphs
 from
$ B_e $ are selected. As operational bases - the minimal
subsets on inclusion of graphs from $ B_o $
specifying types of gluing operations, applicable to
superpositions of graphs selected earlier from elemental bases.
Minimal on inclusion means that a closed class with the same
characteristic
property cannot be obtained by using any of own subsets from
 $ B_o $.
Subsets of graphs from $ B_o $, not satisfying the specified
constraints for the graphs selected from $ B_e $ correspond to
empty cells in the tables. The cells in the table 1 also remain
empty, if graphs selected from $ B_o $ are not isomorphic to
subgraphs of any graph selected from $ B_e $ (cells with subsets
of $B_e$ not containing $K _2$  and subset of $B_o$ containing
only $O_2$).
$\Box$

\bigskip
\centerline{Table 1.The characteristic properties of the connected
graphs}
\bigskip
\begin{tabular}{|p{1.8cm}|p{3cm}|p{3cm}|p{3.2cm}|}
\hline
&&&\\
 $  B_{e}\, \backslash \, B_{o} $ & $ \quad O_1, O_2$ & $\qquad O_2$ & $\qquad O_1$  \\[10pt]
\hline
&&&\\
$O_1,C_1,K_2$ &\centerline{All connected} \centerline{graphs}
 &
 \centerline{Graph
$C_1$ or} \centerline{multigraphs}  \centerline{with $N \leq 2$}
&\centerline{Graphs} \centerline{without cycles}
\centerline{$C_{n}$, $n \geq 2$}
\\
\hline
&&&\\
\centerline{$ C_1,K_2$} & \centerline{Graphs with} \centerline{$M
\geq 1$} &
 \centerline{Graph
$C_1$ or} \centerline{multigraphs } \centerline{with $N = 2$} &
\centerline{Graphs with} \centerline{$M \geq 1$}
\centerline{without cycles} \centerline{$C_{n}$, $n \geq 2$}
\\
\hline
&&&\\
\centerline{$O_1,K_2$} & \centerline{Multigraphs} &
\centerline{Multigraphs} \centerline{with $N \leq 2$} &\centerline
{Trees}
\\
\hline
&&&\\
$\quad K_2$ &\centerline{Multigraphs} \centerline{with $N \geq 2$}
&\centerline{Multigraphs} \centerline{with $N=2$} &Trees with $N
\geq 2$
\\
\hline
&&&\\
\centerline{$ O_1, C_1$}
 & \centerline {---}
& \centerline {---} & \centerline { Graphs with $N=1$}
\\
\hline
&&&\\
$\quad C_1$
 & \centerline {---}
 & \centerline {---}
&  \centerline{Graphs with  $N=1$}
 \centerline{
and $M \geq 1$ }
\\
\hline
\end{tabular}
\bigskip

Using data from tables 1 and 2, we construct for a closed class of
all graphs the diagram of the inclusions of all its closed
subclasses, being basis precomplete in the relevant superclasses.
 Add for completeness the lower
part of the diagram with four trivial closed classes (Fig. 4).

\newpage
\centerline{Table 2.The characteristic properties of the
disconnected graphs}
\bigskip
\begin{tabular}{|p{1.8cm}|p{2.4cm}|p{2.4cm}|p{2.4cm}|p{2.4cm}|}
\hline
&&&&\\
 $  B_{e}\, \backslash \, B_{o} $ & $ \quad O_0, O_1, O_2$ & $\quad O_0, O_2$ & $\quad O_0, O_1$  & $\qquad O_0$\\[5pt]
\hline
\begin{tabular}{p{1.8cm}}
$\hspace {-0.3cm}O_1,C_1,K_2$
 \end{tabular}
& \begin{tabular}{p{2.4cm}} All graphs
\end{tabular}
& \begin{tabular}{p{2.4cm}}
 \hspace {-0.25cm} No graphs \\ \hspace {-0.3cm} with $N=1$\\\hspace {-0.25cm} and
 $M \geq 2$
\end{tabular}
& \begin{tabular}{p{2.4cm}} \hspace {0.15cm} Graphs \\ \hspace
{0.35cm}without
\\
\hspace {0.35cm}
 cycles
\\\hspace {0.2cm}$C_{n}$, $n \geq 2$
\end{tabular}
&\begin{tabular}{p{2.4cm}} \hspace {-0.2cm}Connectivity
\\\hspace {-0.15cm}components \\
\hspace {-0.2cm}isomorphic to\\\hspace {-0.25cm} \hspace
{-0.1cm}$O_1 \lor C_1\lor K_2$
\end{tabular}
\\
\hline
\begin{tabular}{p{1.8cm}}
$ C_1,K_2$
 \end{tabular}
& \begin{tabular}{p{2.4cm}}
\hspace {0.35 cm}Graphs\\ \hspace {0.35 cm}without\\
\hspace {0.25cm}
 isolated\\\hspace {0.25cm} vertices
\end{tabular}
& \begin{tabular}{p{2.4cm}}
 \hspace {-0.25cm}Graphs with \\\hspace {-0.25cm}perfect edge \\
  matching \footnotemark{}

\end{tabular}

& \begin{tabular}{p{2.4cm}}
 \hspace {0.25cm}Graphs \\\hspace {0.25cm}without \\
\hspace {0.25cm}isolated\\ \hspace {0.2cm} vertices \\
and cycles
\\ \hspace {0.1cm}$C_{n}$, $n \geq 2$
\end{tabular}
&\begin{tabular}{p{2.4cm}}

\hspace {-0.2cm}Connectivity
\\\hspace {-0.15cm}components \\
\hspace {-0.2cm}isomorphic to
\\ \hspace {0.1cm}$C_1 \lor K_2$
\end{tabular}
\\
\hline
\begin{tabular}{p{1.8cm}}
$O_1,K_2$
 \end{tabular}
& \centerline {---} & \begin{tabular}{p{2.4cm}} \hspace {-0.2cm}
Multigraphs
\end{tabular}
& \begin{tabular}{p{2.4cm}} \hspace {0.5cm}Woods
\end{tabular}
&\begin{tabular}{p{2.4cm}} \hspace {-0.2cm}Connectivity
\\\hspace {-0.15cm}components \\
\hspace {-0.2cm}isomorphic to \\\hspace {0.15cm} $O_1\lor K_2$
\end{tabular}
\\
\hline
\begin{tabular}{p{1.8cm}}
$K_2$
 \end{tabular}
& \begin{tabular}{p{2.4cm}} \hspace {-0.3cm} Multigraphs\\
\hspace{0.1cm}
 without\\\hspace {0.2cm}izolated\\\hspace {0.08cm}
vertices
\end{tabular}
& \begin{tabular}{p{2.4cm}} \medskip Multigraphs \\
with perfect \\\hspace {-0.25 cm}edge matching
 \end{tabular}
& \begin{tabular}{p{2.4cm}} \hspace {-0.4cm}Woods without
\\\hspace {0.3cm}izolated \\\hspace {0.3cm}vertices
\end{tabular}
&\begin{tabular}{p{2.4cm}} \hspace {-0.2cm}Connectivity
\\\hspace {-0.15cm}components \\
\hspace {-0.2cm}isomorphic to
  \\
 \hspace {0.5cm}$K_2$\\
\end{tabular}
\end{tabular}

\begin{tabular}{|p{1.8cm}|p{2.4cm}|p{2.4cm}|p{2.4cm}|p{2.4cm}|}
\hline &&&&
\\
[-5pt] \hspace {11cm} \centerline{$O_1,C_1$} & \hspace
{2.4cm}\centerline {---} &
 Connectivity
 components with $n=1$ and
  $m
\geq 0 ,$
 if $N=1$
 then $M \leq 1$ &
 Connectivity  \hspace {1.5cm}components with $n=1,$  \centerline{$m
\geq 0$}&

\hspace {-0.1cm}Connectivity  components \hspace {0.2cm}isomorphic
to \hspace {0.1cm}$O_1\lor C_1$
\\
\end{tabular}

\begin{tabular}{|p{1.8cm}|p{2.4cm}|p{2.4cm}|p{2.4cm}|p{2.4cm}|}
\hline
&&&&\\
\hspace {0.4cm} $C_1$ & \centerline {---}
 & Connectivity
   \hspace
{1.5cm}components with $n=1,$ \centerline{$m \geq 1,$}\hspace
{1.5cm} $M-N=2k$,\hspace {1.25cm}$k=0,1,\dots$
 & Connectivity   \hspace {1.5cm}components  with $n=1,$  \centerline{$m \geq
 1$}
 & Connectivity  components  isomorphic to\hspace {0.25cm}
  $C_1$
\\[3pt]
\hline
&&&&\\
\hspace {0.4cm} $O_1$ & \centerline {---} & \centerline {---} &
\centerline {---} & \centerline { \hspace {0.1cm} Empty \hspace
{0.4cm}}\hspace {0.3cm} graphs
\\[5pt]
\hline
\end{tabular}

 \footnotetext{The parity of the sum of the number of loops in
the components
with $ n \geq 2 $ must coincide with the parity of the sum
 $ \sum_ {i=1} ^ {q} (m_ {i} -1) $, where $ q $ is the number
 of components with $ n = 1, \, m_ {i} $ is the number
  of loops in the $i$-th component}

\newpage

    \begin{network}

        \nnNode"1"(4,45)    {+[o][F]{\,1\,}}
        \nnNode"2"(-26,35)    {+[o][F]{\,2\,}}
        \nnNode"3"(-6,35)    {+[o][F]{\,3\,}}
        \nnNode"4"(14,35)    {+[o][F]{\,4\,}}
        \nnNode"5"(34,35)    {+[o][F]{\,5\,}}
        \nnNode"6"(-36.3,20)    {+[o][F]{\,6\,}}
        \nnNode"7"(-29.1,20)    {+[o][F]{\,7\,}}
        \nnNode"8"(-21.8,20)    {+[o][F]{\,8\,}}
        \nnNode"9"(-14.6,20)    {+[o][F]{\,9\,}}
        \nnNode"10"(-7.3,20)    {+[o][F]{10}}
        \nnNode"11"(0,20)    {+[o][F]{11}}
        \nnNode"12"(7.3,20)    {+[o][F]{12}}
        \nnNode"13"(14.6,20)    {+[o][F]{13}}
        \nnNode"14"(21.8,20)    {+[o][F]{14}}
        \nnNode"15"(29.1,20)    {+[o][F]{15}}
        \nnNode"16"(36.3,20)    {+[o][F]{16}}
        \nnNode"17"(43.5,20)    {+[o][F]{17}}
        \nnNode"18"(-36.3,0)    {+[o][F]{18}}
        \nnNode"19"(-29.7,0)    {+[o][F]{19}}
        \nnNode"20"(-23.1,0)    {+[o][F]{20}}
        \nnNode"21"(-16.5,0)    {+[o][F]{21}}
        \nnNode"22"(-9.9,0)    {+[o][F]{22}}
        \nnNode"23"(-3.3,0)    {+[o][F]{23}}
        \nnNode"24"(3.3,0)    {+[o][F]{24}}
        \nnNode"25"(9.9,0)    {+[o][F]{25}}
        \nnNode"26"(16.5,0)    {+[o][F]{26}}
        \nnNode"27"(23.1,0)    {+[o][F]{27}}
        \nnNode"28"(29.7,0)    {+[o][F]{28}}
        \nnNode"29"(36.3,0)    {+[o][F]{29}}
        \nnNode"30"(43.5,0)    {+[o][F]{30}}
        \nnNode"31"(-33,-15)    {+[o][F]{31}}
        \nnNode"32"(-20,-15)    {+[o][F]{32}}
        \nnNode"33"(-7,-15)    {+[o][F]{33}}
        \nnNode"34"(10,-15)    {+[o][F]{34}}
        \nnNode"35"(24,-15)    {+[o][F]{35}}
        \nnNode"36"(37,-15)    {+[o][F]{36}}
        \nnNode"37"(-13,-25)    {+[o][F]{37}}
        \nnNode"38"(4,-25)    {+[o][F]{38}}
        \nnNode"39"(21,-25)    {+[o][F]{39}}
        \nnNode"40"(4,-32)    {+[o][F]{40}}
        \nnLink"1,2"      {@{-}}
        \nnLink"1,3"      {@{-}}
        \nnLink"1,4"      {@{-}}
        \nnLink"1,5"      {@{-}}
        \nnLink"2,6"      {@{-}}
        \nnLink"2,7"      {@{-}}
        \nnLink"2,8"      {@{-}}
        \nnLink"2,9"      {@{-}}
        \nnLink"3,7"      {@{-}}
        \nnLink"3,10"      {@{-}}
        \nnLink"3,11"      {@{-}}
        \nnLink"3,12"      {@{-}}
        \nnLink"4,14"      {@{-}}
        \nnLink"4,13"      {@{-}}
        \nnLink"4,16"      {@{-}}
        \nnLink"4,8"      {@{-}}
        \nnLink"4,10"      {@{-}}
        \nnLink"5,16"      {@{-}}
        \nnLink"5,17"      {@{-}}
        \nnLink"5,11"      {@{-}}
        \nnLink"5,9"      {@{-}}
        \nnLink"6,18"      {@{-}}
        \nnLink"6,22"      {@{-}}
        \nnLink"6,23"      {@{-}}
        \nnLink"7,18"      {@{-}}
        \nnLink"7,19"      {@{-}}
        \nnLink"7,20"      {@{-}}
        \nnLink"8,19"      {@{-}}
        \nnLink"8,21"      {@{-}}
        \nnLink"8,22"      {@{-}}
        \nnLink"9,20"      {@{-}}
        \nnLink"9,23"      {@{-}}
        \nnLink"10,19"      {@{-}}
        \nnLink"10,24"      {@{-}}
        \nnLink"10,27"      {@{-}}
        \nnLink"10,25"      {@{-}}
        \nnLink"11,20"      {@{-}}
        \nnLink"11,24"      {@{-}}
        \nnLink"11,28"      {@{-}}
        \nnLink"12,18"      {@{-}}
        \nnLink"12,25"      {@{-}}
        \nnLink"12,28"      {@{-}}
        \nnLink"14,29"      {@{-}}
        \nnLink"14,21"      {@{-}}
        \nnLink"14,27"      {@{-}}
        \nnLink"13,30"      {@{-}}
        \nnLink"13,25"      {@{-}}
        \nnLink"13,22"      {@{-}}
        \nnLink"16,29"      {@{-}}
        \nnLink"17,28"      {@{-}}
        \nnLink"17,30"      {@{-}}
        \nnLink"17,23"      {@{-}}
        \nnLink"18,31"      {@{-}}
        \nnLink"18,32"      {@{-}}
        \nnLink"19,33"      {@{-}}
        \nnLink"20,32"      {@{-}}
        \nnLink"21,33"      {@{-}}
        \nnLink"22,31"      {@{-}}
        \nnLink"23,32"      {@{-}}
        \nnLink"24,34"      {@{-}}
        \nnLink"24,35"      {@{-}}
        \nnLink"25,31"      {@{-}}
        \nnLink"25,35"      {@{-}}
        \nnLink"27,34"      {@{-}}
        \nnLink"26,33"      {@{-}}
        \nnLink"28,35"      {@{-}}
        \nnLink"28,32"      {@{-}}
        \nnLink"29,36"      {@{-}}
        \nnLink"29,34"      {@{-}}
        \nnLink"30,36"      {@{-}}
        \nnLink"30,35"      {@{-}}
        \nnLink"30,36"      {@{-}}
        \nnLink"30,35"      {@{-}}
        \nnLink"31,37"      {@{-}}
        \nnLink"32,37"      {@{-}}
        \nnLink"33,38"      {@{-}}
        \nnLink"34,38"      {@{-}}
        \nnLink"35,37"      {@{-}}
        \nnLink"36,39"      {@{-}}
        \nnLink"37,40"      {@{-}}
        \nnLink"38,40"      {@{-}}
        \nnLink"39,40"      {@{-}}
        \nnLink"11,26"      {@{-}}
        \nnLink"15,5"       {@{-}}
        \nnLink"29,15"      {@{-}}
        \nnLink"26,15"      {@{-}}
        \nnLink"34,26"      {@{-}}
        \nnLink"31,19"      {@{-}}
        \nnLink"24,16"      {@{-}}
   \end{network}
\smallskip
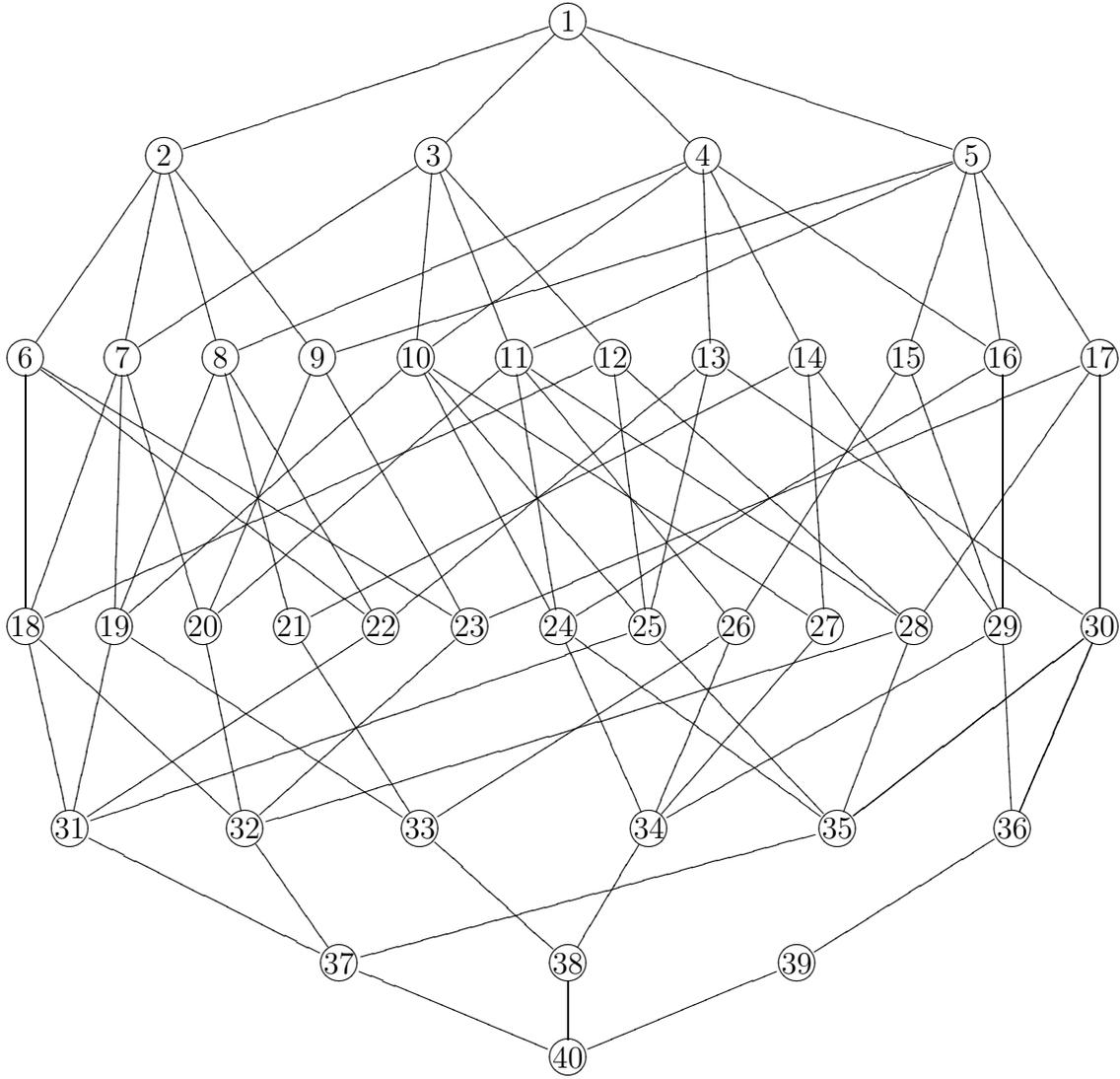
 \centerline{ Figure 4. The diagram of inclusions of all
basis precomplete subclasses of the set of all graphs. }
\medskip

 Denote the subclasses of the set $ \Im$ with
generator bases $B'_e \subseteq B_e $ and $ B'_o \subseteq B_o $
as $ \Im (B'_e, B'_o) $. The generating bases are given below.
\medskip
\\
1. $\Im(\{O_1,C_1,K_2\},\{O_0,O_1,O_2\})$.\,
 2. $\Im(\{O_1,C_1,K_2\}, \{O_1,O_2\})$.\\
 3. $\Im(\{C_1,K_2\}, \{O_0,O_1,O_2\}). \, $4. $\Im(\{O_1,C_1,K_2\},\{O_0,O_1\}).$
5. $\Im(\{O_1,C_1,K_2\},\{O_0,O_2\})$. \\6. $\Im(\{O_1,K_2\},
\{O_1,O_2\})$.\, 7. $\Im(\{C_1,K_2\}, \{O_1,O_2\})$.\, 8.
$\Im(\{O_1,C_1,K_2\}, \{O_1\})$.\\ 9. $\Im(\{O_1,C_1,K_2\},
\{O_2\})$.\, 10. $\Im(\{C_1,K_2\}, \{O_0,O_1\})$.\,
 11. $\Im(\{C_1,K_2\},
\{O_0,O_2\})$.\\ 12. $\Im(\{K_2\}, \{O_0,O_1,O_2\})$.\, 13.
$\Im(\{O_1,K_2\}, \{O_0,O_1\})$. \, 14. $\Im(\{O_1,C_1\},
\{O_0,O_1\})$ \\ 15. $\Im(\{O_1,C_1\},\{O_0,O_2\})$. 16.\,
$\Im(\{O_1,C_1,K_2\}, \{O_0 \})$.\, 17. $\Im(\{O_1,K_2
\},\{O_0,O_2\})$.\\ 18. $\Im(\{K_2\}, \{O_1,O_2\})$.\, 19.\,
$\Im(\{C_1,K_2\}, \{O_1\})$.\,
 20. $\Im(\{C_1,K_2\}, \{O_2\})$.\\ 21. $\Im(\{O_1,C_1\}, \{O_1\})$.\, 22.
$\Im(\{O_1,K_2\}, \{O_1\})$.\,  23. $\Im(\{O_1,K_2\}, \{O_2\})$.\,
\\24. $\Im(\{C_1,K_2\}, \{O_0 \})$.\, 25. $\Im(\{K_2\}, \{O_0,O_1
\})$. \, 26. $\Im(\{C_1\}, \{O_0,O_2 \})$.
\\ 27. $\Im(\{C_1\}, \{O_0,O_1 \})$.\,
 28. $\Im(\{K_2\},
\{O_0,O_2\})$. \, 29. $\Im(\{O_1,C_1\}, \{O_0 \})$. \\ 30.
$\Im(\{O_1,K_2\}, \{O_0 \})$. \, 31. $\Im(\{K_2\}, \{O_1 \})$.\,
 32. $\Im(\{K_2\},
\{O_2\})$. \, 33. $\Im(\{C_1\}, \{O_1 \})$. \\34. $\Im(\{C_1\},
\{O_0 \})$. \, 35. $\Im(\{K_2\}, \{O_0 \})$.\,
 36. $\Im(\{O_1\}, \{O_0\})$.\, 37. $\Im(\{K_2\},\{K_2 \})$. \\38. $\Im(\{C_1\},\{C_1 \})$.
 39. $\Im(\{O_1\},\{O_1 \})$. \, 40. $\Im(\{O_0\},\{O_0 \})$.
 \medskip

The graphs from the considered closed classes possess the most
"strong" \ characteristic properties because for their
constructive descriptions enough using finite elemental and
operational bases.

Constructive descriptions of closed classes with more "weak" \
 characteristic properties also include restrictions on the
choice
of identified subgraphs in the operand graphs and can be on
the
method of identification (situations corresponding to Fig.2
 and
Fig.3). Gluing operations satisfying such restrictions are
denoted as $ H $-gluing operations.  A class of graph
closed
with respect to the operations of $ H $-gluing is called for
brevity an $ H $-closed class.

These restrictions determine the {\it dynamic
characterization} of classes of graphs with a given property.
Together with generating bases, they give constructive
descriptions of closed classes of graphs.
  \medskip

\leftline {{\bf 2.3. Constructive
descriptions  of closed classes of graphs.}}
\medskip

We will consider closed classes of graphs
with some classical properties.
\medskip

\leftline {{\bf 2.3.1. Triangulated graphs. }}
\medskip

A graph $ G $ is called {\it triangulated, or
chordal} if it does not contain a simple cycle
$ C_n, n \geq 4 $ without a chord - edge connecting
non-adjacent vertices of a cycle.

\centerline {{\bf Dynamic characterization }}

From the definition of triangulated graphs it follows that the
presence or absence in the graph of multiple edges does not affect
the triangulation property, then we restrict ourselves to
considering simple triangulated graphs.

\medskip
{\bf Lemma 2.3.1.1} [7]. {\it Operations over the complete
subgraphs of gluing preserve the triangulation of graphs.}

\centerline {{\bf Generating bases}}

The operations preserve the triangulation of graphs are denoted as
the \mbox{$H_t$-gluing}.
\medskip

{\bf Theorem 2.3.1.1} [7]. {\it $H_t$-closed class of triangulated
graphs has countable generating bases $ B_e = \{O_1, K_2, K_3,
\ldots \} \ and \ B_o =\{O_1, K_2, K_3, \ldots \} $}
\medskip

\leftline {{\bf 2.3.2. Planar graphs. }}
\medskip

A graph is called planar if it admits a geometric
implementation on the plane, that is, the vertices of the
graph can
arrange on the plane so that none of its edges intersect
and do not go through extraneous vertices.

\centerline {{\bf Dynamic characterization }}

Two characterizations of planar graphs were considered: based on
traditional geometric representations and using only the
set-theoretic approach based on constructive descripti\-ons of
graphs.

\centerline {Geometric representations}
\medskip

{\bf Lemma 2.3.2.1} [7].
{\it Each planar graph $ G $ can be flat packed into
where all the vertices of an arbitrary face $ f $ with a
connected boundary
are located on a circle inscribed in the face $ f $
in the order of circular traversal of faces.}
\medskip

Suppose that all vertices from $ V (G'_1)$ and $ V (G'_2) $ belong
to flat stacking of planar graphs $ G_1 $ and $ G_2 $ respectively
to faces $ f_1 $ and $ f_2 $ with connected boundaries. Convert
flat styling graphs $ G_1 $ and $ G_2 $ so that all vertices of
the faces $ f_1 $ and $ f_2 $ were located on a circle inscribed
in the face $ f $ in the order of circular traversal of faces.
 We identify the subgraphs $ G'_1$ and $
G'_2$, choosing pairs identified vertices in accordance with
circular rounds of these circles. Gluing operations matching the
specified restrictions on the choice and method of identifying the
subgraphs $ G'_1$ and $ G'_2$ are denoted as
operations of $H_p$-gluing.
\medskip

 {\bf Lemma 2.3.2.2} [7].
{\it $H_p$-gluing operations preserve planar graphs.}
\medskip

\centerline { Set-theoretic approach}

Let $G'\subset G $.
The subgraph of the graph $ G $ generated by the edges of
 set of $
E (G) \backslash E(G')$ is the {\it shell} of subgraph
$ G'$ in
the graph $ G $.
A connected subgraph of a planar graph $ G $, all the vertices
and edges of which belong to the boundary of any connected
face $ f $, is denoted by $ G_f $.

The graph $ G_f $ is represented in the form of superposition
of simple cycles and trees with gluing
subgraphs $O_1$.
The chain connecting two vertices of
the same simple cycle is called {\it
chordal} if its edges and interior vertices do not belong to
the cycle. $ G_f \subset G $ is
the graph of {\it maximal face} if its shell
 in $ G $ consists only from chordal chains.
\medskip

{\bf Lemma 2.3.2.3} [12].
{\it The subgraph $ G'$ of the planar graph
$ G $ is
graph of a maximal face in some
 plane
packing of the connected graph $ G $ if and only if:

1. $ G'$ is realized by a superposition of simple cycles
 and trees with gluing subgraphs $ O_1 $;

2. $ G'$ is selected in $ G $ so that:

a) the shell $ G'$ consists of a set of chordal chains;

b) each pair of chordal chains connecting the
vertices,
arranged along the cycle in alternating order,
has a common inner vertex;

c) no three vertices of one cycle connect to the
 two vertices from the shell of the subgraph
$ G'$ by disjoint chains.}
\medskip

 Consider the depth first search procedure, in which
the following restrictions
are used:

 1) do not select edges that are bridges in
unfulfilled subgraph of the source graph $ G'$, if there are
other possibilities;

2) among the edges that are bridges, do not
choose those that belong
 to the two connected components
of the source
graph $ G'$, if there are other possibilities.

This procedure is called {\it $db$-search}. In [13]
it was shown that the numbering of the vertices of a planar
graph, realized by
superposition of trees and simple cycles with subgraphs of
gluing $ O_1 $,
in accordance with a $db$-search allows single page
flat lay.
\smallskip

 {\bf Theorem 2.3.2.1} [14].
 {\it Each planar graph $ G $ admits a flat packing in which
 all vertices of the face $f$ are located on the circle
 inscribed in the face $f$,
 in accordance with any $db$-search on the subgraph $ G_f$.}
\medskip

\medskip
 {\bf Theorem 2.3.2.2} [1,15].
 {\it
The graph $ G = (G_1 \circ G_2) \tilde G $
preserves the
planarity of the operand graphs
 $ G_1 $ and $ G_2 $ if:

 1) $ G'_1 \subseteq G_1 ^ {f} \subseteq G_1 $ and
$ G'_2 \subseteq G_2 ^ {f} \subseteq G_2,
 G'_1 \cong G'_2 $
and the subgraphs $ G_1 ^ {f} $ and $ G_2 ^ {f} $ have
the following properties:

a)they are realized by superpositions of simple cycles
 and trees with gluing operations at $ O_1 $;

b) their shells consist of sets of chordal chains;

c) each pair of chordal chains of the shell connecting the
 vertices,
arranged along the cycle in alternating order, has a common
inner vertex;

d) no three vertices of one cycle connect to the two
vertices of the shell by disjoint chains.

2. Pairs of identifiable vertices of the subgraphs
 $ G'_1$ and
 $ G'_2$ are selected according to their order
enumerations for arbitrary $db$-search  for $ G_1 ^ {f} $
and $ G_2 ^ {f} $. Identifiable
 pairs
 edges are selected from the set of multiple edges formed
 as a result of the identification of the vertices.}.
\medskip

So, the restrictions on operations of $H_p$-gluing preserving
the
planarity of graphes can be formulated using only
set-theoretic approach.
\medskip

\centerline {\bf {Generating bases }}

The generating bases of the closed class of all planar
graphs coincide with the bases of the class of all graphs.
Only the above restrictions are introduced on the choice of
identifiable subgraphs in the operand graphs and on the
identification method.
\medskip

 {\bf  Theorem 2.3.2.3} [7].
 {\it The $H_p$-closed class of planar graphs has an
 elemental  basis
 $ B_e = \{O_1, C_1, K_2 \}$ and the operational basis
$ B_o = \{O_0, O_1, O_2 \}$.}
\medskip

We restrict ourselves further to the consideration of simple
planar graphs. The gluing operation preserves the absence of
multiple edges,
 if each
a pair of vertices non-adjacent in $ \tilde G $ corresponds
to a pair of non-adjacent
vertices in at least one of the operand graphs $ G_1 $ or
$ G_2 $. Such
gluing operations are referred to as
{\it $ \prec H \succ $-gluing operations}.
\medskip

 {\bf  Corollary 2.3.2.1.}
 {\it
 $ \prec H_p \succ $-closed  class of simple planar
 graphs has elemen\-tal
 basis $ B_e = \{O_1, K_2 \} $ and operational basis
 $ B_o = \{O_0, O_2 \} $.}
\medskip

All restrictions on $H$-gluing operation can be divided into
internal ones, without
which gluing operations cannot preserve the required
characteristic
property of graphs (note that all previously used
restrictions  were
internal ones), and
additional external ones that affect the
power of generating bases, order of graph assembly, the
amount of
redundancy of the constructive description, etc.

Consider as an example
gluing operations in which the set vertices of the subgraph
gluing $ V (\tilde G)$
is separating set in the resulting graph.
 Denote such operations as
$H_s$-gluing. Then, for the class of simple planar graphs we
have
 $ \prec H_{ps} \succ $-gluing operations.
\medskip

 {\bf Theorem 2.3.2.4} [16].
 {\it
$ \prec H_{ps} \succ $-closed class of simple planar graphs has
the generat\-ing bases $ B_e = \{O_1, K_2, K_3, K_4 \} $ and
 $ B_o = \{O_0, O_1,
O_2, O_3, O_4, O_5 \} $.}
 \smallskip

If we go to gluing by the generated subgraphs ($<H>$-gluing) and
require that the sets of vertices $ V(\tilde G)$ are minimal
separating in the resulting graphs then, for the class of simple
planar graphs, $<H_{pv}>$-gluing operations should be used.
 \medskip

{\bf Theorem 2.3.2.5} [17].
 {\it
The $<H_{pv}>$-closed class of simple planar graphs has
elemental
basis
 \mbox {$ B_e = \{O_1, K_2, K_3, K_4 \} $.}
 Operational basis $ B_o $
contains  \hbox {16} types operations, whose gluing subgraphs are
isomorphic to graphs from the set }\\ [- 5mm]
 \bigskip
\[
\begin{array}{c}
\{\ O_0,\, O_1,\, O_2,\, K_2,\, O_3,\, (O_1\circ K_2)O_0,\,
 L_3,\, K_3,\, O_4,\smallskip \\[-1mm]
(O_2\circ K_2)O_0,\, (K_2\circ K_2)O_0,\,
 (O_1\circ L_3)O_0,\,
L_4,\, C_4,\,\, L_5,\, C_5 \ \}. \\
\end{array}
\]
\medskip
The operational basis will increase even more if you restrict
using of gluing operations so that the number of edges in
$E(\tilde G)$ is also minimal. The correspond\-ing
 operations denoted
as operations of $<H_{pve}>$-gluing.
\medskip

{\bf Theorem 2.3.2.6} [7,18].
 {\it
The $<H_{pve}>$-closed class of simple planar graphs has an
elemental basis
 \mbox {$ B_e = \{O_1, K_2, K_3, K_4 \} $.}
Operational basis $ B_o $ contains \hbox {22} types operati\-ons,
whose gluing subgraphs are isomorphic to graphs from the set }\\
[- 5mm]
\smallskip
\[
\begin {array} {c}
\{\ O_0, \, O_1, \, O_2, \, K_2, \, O_3, \, (O_1 \circ K_2)
 O_0, \,
 L_3, \, K_3, \, O_4, \smallskip\\ [- 1mm]
(O_2 \circ K_2) O_0, \, (K_2 \circ K_2) O_0, \,
(O_1 \circ L_3) O_0, \,
L_4, \, C_4, \, O_5, \, (O_3 \circ K_2) O_0,\smallskip \\ [- 1mm]
((K_2 \circ K_2) O_0 \circ O_1) O_0, \, (O_2 \circ L_3) O_0, \,
 (O_1 \circ L_4) O_0, \, (K_2 \circ L_3) O_0, \, L_5, \, C_5
\ \}. \\
\end {array}
\]
\medskip

\centerline {\bf{Triangulated planar graphs}}
\medskip

{\bf Theorem 2.3.2.7} [7]. {\it $<H_{p}>$-closed class of
triangulated simple planar graphs has such bases
 \mbox {$ B_e = \{O_1, K_2, K_3, K_4 \} $} and
 $ B_o = \{O_0, O_1, K_2, K_3 \} $.}
\medskip

\centerline {\bf{Maximality planar graphs}}
\medskip

Simple planar graph $ G $ is called {\it maximal}, if adding
any edge to $ G $ takes it out of class planar.
\medskip

 {\bf  Theorem 2.3.2.8} [7].
{\it
$<H_{p}>$-closed class of maximal planar graphs has
countable elemental basis $ B_e $ and operational basis
 $ B_o = \{K_3 \} $.}
\medskip

Constructive descriptions of closed classes of outerplanar as well
as triangula\-ted and maximal outerplanar graphs are given in
[19].
\medskip

\leftline {{\bf 2.3.3. Euler graphs. }}
\medskip

A connected graph $ G $ with even degrees of vertices is
 called  {\it Euler}.

\centerline {{\bf Dynamic characterization }}
\medskip

{\bf Lemma 2.3.3.1} [20]. {\it If the graphs $ G_1 $ and $ G_2 $
are Euler, then the resulting
 graph $ G = (G_1 \circ G_2) \tilde G $ will be Euler if and
only if the degrees
all vertices of the gluing subgraph $ \tilde G $ are even.}
\medskip

In order to reduce the redundancy of constructive descriptions
of Euler
graphs we restrict ourselves to the use of glue operations on
empty
subgraphs. We denote them as operations of the
$ H^{\emptyset} $-gluing.

\centerline {{\bf Generating bases }}
\medskip

 {\bf Theorem 2.3.3.1} [21].
{\it The $H^{\emptyset} $-closed class of Euler graphs has
generating bases $ B_e = \{C_1, C_2, ... \} $ and
$ B_o = \{O_1, O_2,... \} $.}
\medskip

The infinity of the operational basis follows from the result
 of Alon [22].

\medskip

\leftline {{\bf 2.3.4. Euler planar graphs.}}
\medskip

Operations of gluing preserve the euler and planar properties of
graphs are denoted as $H^{\emptyset}_p$-gluing. Based on studies
performed in [8,20,21,23,24], the following result was obtained.

\medskip

 {\bf Theorem 2.3.4.1.}
 {\it
$H^{\emptyset}_p$-closed class of Euler planar graphs has the
elemental basis \linebreak $ B_e = \{C_1,C_2, ...\}$ and three
operating bases $B^{1}_{o} = \{O_1, O_2, O_3 \} $, \mbox{$ B^{2}_o
= \{O_1, O_2, O_4 \} $} and $ B^{3}_o = \{O_1, O_2, O_5 \} $.}
\medskip

Graphs that do not contain vertices of the second degree are
called  {\it topological}. If the gluing operations are
carried
out on the generated subgraphs, then the operational basis of
the
closed class of simple topological Euler planar graphs
becomes infinite [25].
\medskip

\leftline {{\bf 2.3.5. Hamiltonian graphs.}}
\medskip

A graph $G$ is called {\it Hamiltonian} if it is possible to
select in it a cycle containing all the vertices of the graph.

\centerline {{\bf Dynamic characterization }}
\medskip

{\bf Lemma 2.3.5.1} [26]. {\it If $ G_ {1} $ and $ G_ {2} $
are
Hamiltonian graphs, then the resulting graph
$G = (G_ {1} \, \circ \, G_ {2}) \, \tilde {G}$ will also be
Hamiltonian under any of the following conditions:

 1) the identified subgraph of at least one of the operand
graphs contains all its vertices;

 2) the identifiable subgraphs of the operand graphs consist
 of two vertices that are adjacent in their Hamiltonian cycles.}
\medskip

Gluing operations satisfying any of these restrictions
are called $ H_g $-gluing operations.
Since the presence or absence in the graph of multiple edges
 does not affect the Hamiltoni\-an property, then
we restrict ourselves to using operations
 $ \prec H_g \succ$-glues
excluding the appearance of multiple edges.

\centerline{{\bf Generating bases }}

Here, as well as for Euler planar graphs, there are three
operational bases.
\medskip

 {\bf Theorem 2.3.5.1} [26].
{\it The $ \prec H_g \succ $-closed class of Hamiltonian graphs
has elemental basis $ B_ {e} = \{C_1, C_2, \ldots \}$ and three
following operational bases \mbox{$ B_{o} ^ 1 = \{O_1, K_2, C_4,
C_5
 \ldots \} $},
$B_ {o} ^ 2 = \{O_1, K_2, L_3, L_4, \ldots \}
 $
and \,
 $ B_ {o} ^ 3 = \{O_1, K_2, (L_{n'} \, \circ \,
L_{n''}) \,
 O_0 \}, n', n'' \geq 2 $.}
\medskip

If canonical superpositions are admissible when constructing
graphs of some $H$-closed class, then such a class is briefly
called as canonical $H$-closed class.

 {\bf Corollary 2.3.5.1}
{\it Class of Hamiltonian graphs canonically
$ \prec H_g \succ $-closed  with elemental basis $ B_ {e} = \{C_1,
C_2, \ldots \} $ and two operational bases $ B_ {o} ^ 2 = \!
\{O_1, \! K_2, \! L_3, \! L_4, \ldots \! \} $
 \! or
$ B_ {o} ^ 3 = \{O_1, K_2, (L_ {n'} \,
 \circ \, L_ {n''}) \, O_0 \}, n', n'' \geq 2 $.}
\medskip

\leftline {{\bf 2.3.6. Bipartite graphs.}}
\medskip

A graph $ G $ is called {\it bipartite} if
there exists a partition of the set of its
vertices $ V (G) $ into two subsets $ V_ {1} $
and $ V_ {2} $, each of which generates an
empty graph. If the graph $ G $ is not empty,
then the ends each edge $ e \in E (G) $
belongs to different parts.

\centerline {{\bf Dynamic characterization }}

\medskip
{ \bf Lemma 2.3.6.1} [27]. {\it If $ G_{1} $ and $ G_{2} $ are
bipartite graphs, then the result of gluing them together
graph $G
= (G_ {1} \, \circ \, G_ {2}) \, \tilde {G}, | V (\tilde {G})
|\geq 2 $ is bipartite if and only if when for any vertices
$ v_1,
v_2 \in
 V (\tilde {G}) $, connected by chains in
$ G (E_1 \backslash \tilde {E}) $ and
 $ G (E_2
\backslash \tilde {E}) $, lengths these chains have the same
parity.}
\medskip

This restriction on gluing operations is denoted by $ H_{b} $.
\medskip
\newpage
\centerline {{\bf Generating bases }}
\medskip

{\bf Theorem 2.3.6.1} [27]. {\it The class of bipartite graphs $
H_{b}$-closed with an elemental basis
 $ B_ {e} = \{O_1, K_2 \} $ and the operational basis
 $ B_ {o} = \{O_0, O_2 \} $.}

\medskip

  \leftline{{\bf{\large 3. Conclusion}}}
\medskip
\medskip

Concluding the review, we note the following points.

1.The constructive descriptions of graphs show the efficiency of
using gluing operations to uniformly formulate the conditions for
preservation of various characteristic properties of graphs in
terms of restrictions on the type of identified subgraphs, their
choice in operand graphs and the identification method.

The "payment" \, for this universality is the redundancy
introduced by gluing operations in the information about the graph
with labelled vertices. Estimates of the magnitude of this
redundancy are obtained for the Eulerian graphs [28], some classes
of triangulated graphs [29] and Hamiltonian graphs [30,31].

When considering unlabelled graphs, knowledge of their
construction processes
 can significantly reduce the length of
the graph code and complexity of decoding algorithms
 by using the numbering of the vertices,
reflecting graph assembly order [32,33].

2. The presence of several operational bases for some closed
classes of graphs allows you to formulate tasks
of optimal graph synthesis. For example, in the
classical statement of minimizing the number of gluing
operations needed to build a graph.

Another class of optimal graph synthesis problems arises in
supercomputer physical-mathematical modeling design of
large graphs at the stage when it is necessary consistent
 return
from the reduced graph to the original graph of large
dimension with preserving the solution obtained on the
previous steps.

This process can be implemented using subgraphs
duplication operations with full or partial preserving
their neighborhoods in the current graph.
Such operations named as {\it cloning} operations
are discussed in [34,35].

The task of optimal graph synthesis is put here as follows:
based
on the
graph a small dimension it is necessary to construct a graph
 of
large dimension with specified properties for minimum number
of cloning operations. Such a
task considered for trees and bipartite graphs in [36].

3. Constructive approach methodology to solving applied tasks
on graphs will be successful if implemented the following
principles:

- choose as the source graphs a complete system
of graphs (not necessarily elemental basis)
for each of them the considered problem is solved most
effectively

- choose restrictions on the gluing operations so that the
 resulting graph can be built using canonical superposition
("brick by brick").  Its simplified
structure analysis of graphs and therefore, finding the
solution.

Using the above methodology illustrated in [1] on examples of
solving applied problems economical coding and optimal linear
placement of graphs.

\medskip
\medskip
  \leftline{{\bf{\large Acknowledgement}}}
\medskip
\medskip

This research did not receive any specific grant from funding
agencies in the public, commercial, or not-for-profit sectors.
The
preparation of the monograph [1] was supported by the Russian
Foundation for Basic Research (grant number 16-01-00117).
\newpage
\medskip

\textbf{References}
\medskip

[1] Iordanski \; M. \, A. Constructive graph theory and its
applications.  Nizhny Novgorod: Publishing house "Cyrillic",
 2016.
172\,p.(Russian)

[2] Iordanski  M.  A.  About  operations on graphs
// Discrete models in the theory of control
systems //~ Proceedings of the IV International Conference,
 Krasnovidovo, June 19-25, 2000. / Edited by V.B.
 Alekseev, B. A. Zakharova. "--- M .:
     MAX Press, 2000. "--- P. \, ~33--34.(Russian)

[3] Iordanski M. A. Constructive descriptions of graphs
//~ Proceedings
of the XI International Seminar-school "Synthesis and
Complexity of Control Systems"(Nizhny Novgorod, November
 20-25, 2000). Part I. "--- M.: Publishing House of the center
for applied research on Department of the mechanical-mathematical
Faculty of Moscow State University, 2001."--- P.~80--84.(Russian)

[4]  Iordanski M. A. The structure and methods of generating
closed classes of graphs //~ Proceedings of the VII International
Seminar "Discrete Mathematics and Its applicati\-ons" (Moscow,
January 29 - February 2, 2001). Part II. / Edited by O.B.
Lupanov."--- M.: Publishing House of the center for applied
research on Department of the mechanical-mathematical Faculty of
Moscow State University, 2001. "---P.~218--221.(Russian)

[5] Iordanski M. A. The structure and methods of generating
closed classes of graphs// Discrete Mathematics.
 2003. Vol. 15, Iss.  3. "---P.~105--116.(Russian)

[6] Iordanski M. A. Functional approach to the representation
of graphs ~ // Reports of the Russian Academy of Sciences. "---
1997. "--- T. \, 353, \ No \, 3. " --- P. 303-305.(Russian)

[7] Iordanski M.  A. Constructive graph
descriptions // Discrete analysis and operation research.
 1996.
 Vol. 3, Iss.  4. "---P.~35--63.(Russian)

[8] Burkov E. V.  Operational bases of closed classes of graphs~//
Proceedings of the IX International Seminar "Discrete Mathematics
and Its applications", Moscow, June 18-23, 2007. "---M.:
Publishing House of the Moscow State University M.V. Lomonosov,
2007. "---P.~105--116.(Russian)

[9] Iordanski M. A. Constructive classification graphs
 // Modeling and analysis of information systems.
      2012.  Vol.19,  Iss.  4. "---P.~144--153.(Russian)

[10] Iordanski M. A. Some questions of analysis and synthesis
of graphs
 //~ Proceedings of the First International
Conference "Mathematical Algorithms" (Nizhny Novgorod, August
16-19, 1994). / Edited by M.A. Antonets, V. E.
 Alekseev, V.N. Shevchenko. "--- Nizhny Novgorod: Publishing
 House of the Nizhny Novgorod University, 1995.
 "--- P. \, ~33--38.(Russian)

[11] Yablonsky S.V. Introduction to discrete mathematics.
  "---M.: Publishing House "Science", 2001. 384\,p.(Russian)

[12] Iordanski M. A. Set-theoretic description of subgraphs of
faces of planar graphs // Problems of Theoretical cybernetics.
Materials of the XIV International Conference (Penza, May 23-28,
2005)."--- M.: Publishing House of the Department of the
mechanical-mathematical Faculty of Moscow State University,
2005."--- P.~56.(Russian)

[13] Iordanski M. A. Depth Search and single page graph stacking
// Problems of Theoretical cybernetics. Materials of the XIII
International Conference (Kazan, May 27-31, 2002). Part I / Edited
by O.B. Lupanov. "--- M.: Publishing House of the center for
applied research on Department of the mechanical-mathematical
Faculty of Moscow State University, 2002. "---P.~76.(Russian)

[14] Iordanski M. A. Properties of flat stackings of planar graphs
//~ Proceedings of the VII International Seminar "Discrete
Mathematics and Its applications" (Pokrovskoye, March 4-6,
2006)."--- M.: MAX Press, 2006."---
 P.~136--138.(Russian)

[15] Iordanski M. A. Inductive description of the class of
 planar graphs //~ Proceedings
of the XVI International Seminar-school "Synthesis and Complexity
of Control Systems" \\(St. Petersburg, June 26-30, 2006)."--- M.:
Publishing House of the Department of the mechanical-mathematical
Faculty of Moscow State University, 2006. "---P.~41--43.(Russian)

[16] Iordanski M. A. The complexity of constructive
descriptions of planar graphs
//~ Proceedings
of the IX International Seminar-school "Synthesis and Complexity
of Control Systems"(Nizhny Novgorod, December 16-19, 1998)."---
M.: Publishing House of the Department of the
mechanical-mathematical Faculty of Moscow State University,
1999."--- P.~20--24.(Russian)

[17] Iordanski  M. A. Bases of planar graphs~//
 Discrete models in the theory of control systems:
    V International Conference: (Ratmino, May 26-29, 2003).
 "---M.: Publishing Depart-\linebreak
  ment of the Faculty of Computational
Mathematics and Cybernetics Moscow State University M.V.
Lomonosov, 2003."--- P. 36--38.(Russian)

[18] Iordanski M. A. Constructive descriptions of planar
graphs // Problems of Theoretical
cybernetics. Abstracts of the XI International Conference
(Ulyanovsk, June 10-14, 1996)."--- M.: Publishing House of the
Russian State University for the Humanities, 1996.
 "---P.~76--77.(Russian)

[19] Iordanski M. A. Algorithmic descriptions of outerplanar
 graphs
 //~ Proceedings of the Second International
Conference "Matematical Algorithms" (Nizhny Novgorod,
June 26 - July 1, 1995). / Edited by M.A. Antonets, V. E.
 Alekseev, V.N. Shevchenko. "--- Nizhny Novgorod: Publishing
 House of the Nizhny Novgorod University, 1997.
 "--- P. \, ~78--82.(Russian)

[20] Iordanski  M.  A., Burkov E.  V. Constructive descriptions of
Eulerian planar graphs // Discrete models in the theory of control
systems: VI International Conference: Moscow (December 7-11, 2004)
"---M.: Publishing Department of the Faculty of Computa\\tional
Mathematics and Cybernetics Moscow State University M.V.
Lomonosov. 2004. \linebreak "---P.~167--169.(Russian)

[21] Burkov E. V. Constructive descriptions of planar and Eulerian
 Counts~// Bulletin of Nizhny Novgorod State
 university N.I. Lobachevsky. Mathematics.  2010.  Iss. 5 (1).
 \linebreak  P.~165--~170.(Russian)

[22] Alon N. Tough Ramsey Graphs Without Short Cycles~//Jornal of
Algebraic Combinatorics.
      1995. Vol. 4, Iss. 3."---P.~189--195.

[23] Burkov E. V. Another operational basis for the class of Euler
 planar
 graphs // Problems of Theoretical cybernetics. Abstracts of the XV
International Conference (Kazan, June 2-7, 2008). " --- Kazan:
Publishing house "Fatherland". "--- 2008. \\" --- P. \,
13.(Russian)

[24] Burkov E. V. Short cycles in planar graphs with a minimum
 degree of four~ //
 Bulletin of Nizhny Novgorod State University. Math modeling
     and optimal control. "--- 2009.
   "--- \, Iss. \, 4." --- P. 146--148.(Russian)

[25] Iordanski\; M. \, A.Countable operational basis of
 topological Euler planar graphs // Discrete models in the
theory of control systems: VIII International Conference: Moscow
(april 6-9, 2009): Proceedings."---M.: MAX Press, 2009. "--- P. \,
~127--129.(Russian)

[26] Iordanski M. A. Constructive  descriptions of Hamiltonian
graphs ~// Bulletin of Nizhny Novgorod State
university N.I. Lobachevsky. Mathematics.  2012.  Iss.  3 (1).
 \linebreak  "---P.~137--~140.(Russian)

[27] Iordanski\; M. \, A. Constructive descriptions of bipartite
graphs ~ // Problems of Theoretical cybernetics. Abstracts of the
XV International Conference (Kazan, June 2-7, 2008). - Kazan:
Publishing house "Fatherland". "--- 2008." --- P. \, 44.(Russian)

[28]  Iordanski \; M. \, A. Redundancy of constructive
descriptions of Eulerian graphs // Problems of Theoretical
cybernetics. Materials of the XVII International Conference
(Kazan, June 16-20, 2014). Edited by Yu.I. Zhuravlev.
 "--- Kazan: Fatherland, 2014. \linebreak" --- P. 115-116.
(Russian)

[29] Iordanski \; M. \, A. Redundancy of constructive descriptions
 of (r,s)-trees // Discrete models in the theory of control
systems: IX International Conference, Moscow and Moscow
Region, May 20-22, 2015: Proceedings / Edited by V.B.
 Alekseev, D. S. Romanov, B. R. Danilov. "--- M .:
     MAX Press, 2015. "--- P. \, 90--91.(Russian)

[30]  Iordanski M. A. Redundancy of constructive descriptions of
Hamiltonian graphs //~ Proceedings of the XII International
Seminar "Discrete Mathematics and Its applicati\-ons"\, named
after academician O.B. Lupanov, Moscow, June 20-25, 2016."--- M.:
Publishing House of the Department of the mechanical-mathematical
Faculty of Moscow State University, 2016.
"---P.~290--293.(Russian)

[31]  Iordanski M. A. Redundancy of constructive descriptions of
Hamiltonian planar graphs //~ Proceedings of the XI International
Seminar "Discrete Mathematics and Its applications" (Moscow, June
18-22, 2012)."--- M.: Publishing House of the Department of the
mechanical-mathematical Faculty of Moscow State University,
2012."---P.~285--288.(Russian)

[32] Iordanski  M. A. Constructive descriptions and economical
coding of graphs ~ //Problems of Theoretical cybernetics.
Materials of the XII International Conference (Nizhny Novgorod,
May 17-22, 1999). "--- M.: Publishing House of the Department of
the mechanical-mathematical Faculty of Moscow State University,
1999. "---P.~87.(Russian)

[33] Iordanski \; M. \, A. Constructive descriptions and economical
coding of graphs ~ //
     Bulletin of Nizhny Novgorod State University. Math modeling
     and optimal control. "--- 2000.
   "--- \, Iss. \, 1 (22)." --- P. 88--93.(Russian)

[34]  Iordanski M. A. Cloning graphs  // Problems of Theoretical
cybernetics. Abstracts of the XVIII International Conference
(Penza, June 19-23, 2017)."--- M: MAX Press, 2017."---
P.~108--110.(Russian)

[35]  Iordanski M. A. On a class of graph transformations
// Discrete models in the theory of control
systems: X International Conference: Moscow and Moscow Region
 (May 23-25, 2018) Proceedings / Edited by V.B.
 Alekseev, D. S. Romanov, B. R. Danilov.\linebreak
 "--- M.: MAX Press,
 2018. "---P.~139--142.(Russian)

[36]  Iordanski M. A. On the complexity of graph  synthesis
by
cloning  operations //~ Proceedings of the XIII International
Seminar "Discrete Mathematics and Its
 applica- \linebreak tions"
\, named after academician O.B. Lupanov (Moscow, June 17-22, 2019)
"---M.: Publi\-shing Department of the Faculty of Mechanics and
Mathematics Moscow State University M.V. Lomonosov. 2019. "---P.
220--223.(Russian)

\end{document}